\RequirePackage{fix-cm}
\documentclass[smallextended,envcountsect,final=]{svjour3} 
\smartqed 
\usepackage{graphicx}
\usepackage{amssymb}
\usepackage{amsmath}
\usepackage{algorithm}
\usepackage{algpseudocode}	
\usepackage{algorithmicx}
\usepackage{mathrsfs}
\usepackage{bm}
\usepackage{cite}
\usepackage{xcolor}
\usepackage{pgfplots}
\usepackage{subfigure}
\usepackage{booktabs}
\usepackage{multirow}
\usepackage{url}
\usepackage{graphicx} 
\usepackage{setspace}
\journalname{Noname}
\newcommand\norm[1]{\left\lVert#1\right\rVert}
\newcommand{\Hil}{\mathcal{H}}

\newcommand{\be}{\begin{equation}}
\newcommand{\ee}{\end{equation}}
\newcommand{\bee}{\begin{eqnarray}}
\newcommand{\eee}{\end{eqnarray}}
\newcommand{\bse}{\begin{subequations}}
	\newcommand{\ese}{\end{subequations}}
\newcommand*\diff{\mathop{}\!\mathrm{d}}
\newcommand{\sr}{r_{(t)}}
\newcommand{\SR}{r_{(T)}}
\newcommand{\EGO}{Efficient Global Optimization}
\newcommand{\Ego}{Efficient global optimization}
\newcommand{\ego}{efficient global optimization}
\newcommand{\rX}{(\mathcal{X})}
\makeatletter
\newcommand*{\rom}[1]{\expandafter\@slowromancap\romannumeral #1@}
\makeatother

\spnewtheorem{assump}{Assumption}[section]{\bf}{\it}
\spnewtheorem{conject}{Conjecture}[section]{\bf}{\it}

\usepackage{accents}
\usepackage{amsfonts}       
\usepackage{thmtools,thm-restate}
\usepackage{mathtools}
\newcommand{\gss}{G}

\usepackage{accents}

\allowdisplaybreaks[4]

\newcommand{\ubar}[1]{\underaccent{\bar}{#1}}

\ifodd 1
    \newcommand\rev[1]{{\color{black}#1}}
    \newcommand{\com}[1]{\textbf{\color{red} (COMMENT: #1)}} 
\else
    \newcommand\rev[1]{{#1}}
    \newcommand{\com}[1]{}
\fi

\newcommand{\Proof}{\noindent {\it Proof }$\;\;$}
\begin{document}

\title{Lower Bounds on the Worst-Case Complexity of \EGO}


\author{Wenjie Xu \and Yuning Jiang$^*$ \and Emilio T. Maddalena  \and Colin N. Jones\thanks{$^*$Corresponding author.\\ \email{yuning.jiang@ieee.org}
} }

\institute{
W. Xu, Y. Jiang, E. Maddalena, C. Jones \at
Automatic Control Laboratory\\
EPFL, Switzerland
}

\date{Received: date / Accepted: date}

\maketitle

\begin{abstract}
\Ego{ is} a widely used method for optimizing expensive black-box functions such as tuning hyperparameter, and designing new material, etc. Despite its popularity, less attention has been paid to analyzing the inherent hardness of the problem although, given its extensive use, it is important to understand the fundamental limits of \ego{ algorithms}. In this paper, we study the worst-case complexity of the \ego{ problem} and, in contrast to existing kernel-specific results, we derive a \rev{unified} lower bound for the complexity of \ego{ in} terms of the metric entropy of a ball in its corresponding reproducing kernel Hilbert space~(RKHS). Specifically, we show that if there exists a deterministic algorithm that achieves suboptimality gap smaller than $\epsilon$ for any function $f\in S$ in $T$ function evaluations, 
it is necessary that $T$ is at least $\Omega\left(\frac{\log\mathcal{N}(S(\mathcal{X}), 4\epsilon,\|\cdot\|_\infty)}{\log(\frac{R}{\epsilon})}\right)$, where $\mathcal{N}(\cdot,\cdot,\cdot)$ is the covering number, \rev{$S$ is the ball centered at $0$ with radius $R$ in the RKHS and $S(\mathcal{X})$ is the restriction of $S$ over the feasible set $\mathcal{X}$}. Moreover, we show that this lower bound nearly matches the upper bound attained by non-adaptive search algorithms for the commonly used squared exponential kernel and the Mat\'ern kernel with a large smoothness parameter $\nu$, up to a replacement of $d/2$ by $d$ and a logarithmic term $\log\frac{R}{\epsilon}$. That is to say, our lower bound is nearly optimal for these kernels.  
\end{abstract}


\section{Introduction}
Black-box optimization by sequentially evaluating different candidate solutions without access to gradient information is a pervasive problem. For example, tuning the hyperparameters of machine learning models~\cite{bergstra2012random,snoek2015scalable}, optimizing control system performance~\cite{bansal2017goal,xu2021vabo} and discovering drugs or designing materials~\cite{negoescu2011knowledge,frazier2016bayesian}, etc., can all be formulated as a black-box optimization problem without explicit gradient information. Therefore, \ego~\cite{jones1998efficient,shahriari2015taking}, as a sample-efficient method to solve the expensive black-box optimization problem without explicit gradient information, has recently been receiving much attention. In many applications, e.g., tuning the hyperparameters of a deep neural network, each sample can take significant resources such as time and computation. For such problems, understanding the sample complexity of \ego{ is} of great theoretical interest and practical relevance. 

\bigskip
\noindent
There is a large body of literature on the convergence rates of particular \ego{ algorithms}~\cite{srinivas2012information,vazquez2010convergence,wang2014theoretical,russo2016information,vakili2021optimal,de2012exponential}. Two typical analysis set-ups are the Bayesian and non-Bayesian settings\footnote{The Bayesian setting is typically referred to as Bayesian optimization.}. In the Bayesian setting, the black-box function is assumed to be sampled from a Gaussian process, whereas in the non-Bayesian setting, the black-box function is assumed to be \rev{regular} in the sense of having a \rev{bounded} norm in the corresponding reproducing kernel Hilbert space. 

\bigskip
\noindent
As a complement to convergence analysis of different algorithms, complexity analysis tries to understand the inherent hardness of a problem. Specifically, we are interested in answering the question: \emph{for a class of optimization problems, how many queries to an {oracle}, which returns some information about the function, are necessary to guarantee the identification of a solution with objective value at most $\epsilon$ worse than the optimal value}~\cite{nemirovskij1983problem}? Without a complexity analysis, we cannot tell whether existing algorithms can be improved further in terms of convergence rate. This problem is well studied for convex optimization (e.g., in~\cite{nemirovskij1983problem}), but less well understood for \ego. Existing results on lower bounds are usually kernel-specific~\cite{bull2011convergence,scarlett2017lower,scarlett2018tight,ray2019bayesian,cai2021lower} and cannot be directly applied to general kernels.  

\bigskip
\noindent
Intuitively, the complexity of \ego{ largely} depends on the richness or complexity of the functions inside the corresponding reproducing kernel Hilbert space~(RKHS). As an extreme example, if we adopt the linear kernel, after a finite number of noiseless function evaluations, we can uniquely determine the ground-truth function and hence the optimal solution. To measure the complexity of a set of functions, \emph{metric entropy}~\cite{kolmogorov1959varepsilon} is widely used in learning theory. However, as far as we know, the explicit connection between a complexity measure such as metric entropy for a function set and the problem complexity of \ego{ has} not been established. 

\bigskip
\noindent
This paper focuses on the complexity analysis of \ego{ with} general kernel functions in the non-Bayesian \rev{and noiseless} setting. Although noisy setting is more realistic from practical point of view, it is also critical to consider the noiseless setting from complexity theoretic point of view. The rationale is that the noise may introduce additional statistical complexity to the problem and corrupts the inherent complexity analysis of the \ego. In addition, noiseless setting is not a simple extension of noisy setting, existing analysis under noisy setting~(e.g.,~\cite{cai2021lower,ray2019bayesian,scarlett2017lower,scarlett2018tight}) typically relies on strictly positive noise variance. Simply setting noise variance to zero makes the analysis and results diminish. For example, the noisy bound for Squared Exponential~(SE) kernel in~\cite{scarlett2017lower} is $\Omega(\frac{\sigma^2}{\epsilon^2}\left(\log\frac{B}{\epsilon})^{\frac{d}{2}}\right)$, which is dominated by $\frac{\sigma^2}{\epsilon^2}$, where $\sigma^2$ is the noise variance and $B$ is the function norm upper bound. Simply setting $\sigma=0$ gives a meaningless $\Omega(0)$ bound. Without the analysis under noiseless setting, it is unclear whether this dominant $\frac{\sigma^2}{\epsilon^2}$ term is due to noise or due to the inherent complexity of the RKHS.
\vspace{-0.5cm}
\begin{table}[htbp!]
    \centering
    \scriptsize
    \begin{tabular}{|c|c|c|c|c|}
    \hline 
         \textbf{Works} & \textbf{Noise} &\textbf{SE kernel} & \textbf{Mat\'ern kernel} & \textbf{General kernel} \\
         \hline 
         \cite{bull2011convergence}& No & N/A & $\Omega\left((\frac{1}{\epsilon})^{\frac{d}{\nu}}\right)$ &N/A \\
        \hline 
        \cite{scarlett2017lower}& Yes& $\Omega\left(\frac{\sigma^{2}}{\epsilon^{2}}\left(\log \frac{R}{\epsilon}\right)^{d / 2}\right)$ & $\Omega\left(\frac{\sigma^{2}}{\epsilon^{2}}\left(\frac{R}{\epsilon}\right)^{d / \nu}\right)$ &N/A \\
        \hline 
         Ours& No& $\Omega\left(\left(\log \frac{R}{\epsilon}\right)^{d / 2-1}\right)$ & $\Omega\left(\frac{\left(\frac{R}{\epsilon}\right)^{\frac{d}{\nu+d / 2}}}{\log \frac{R}{\epsilon}}\right)$ & $\Omega\left(\frac{\log \mathcal{N}\left(S(\mathcal{X}), 4 \epsilon,\|\cdot\|_{\infty}\right)}{\log \left(\frac{R}{\epsilon}\right)}\right)$ \\
         \hline 
    \end{tabular}
    \caption{\rev{A summary of the state-of-the-art complexity result for \ego. $\sigma^2$ is the noise variance. $R$ is the function norm upper bound. $d$ is the dimension of input space. $\nu$ is the smoothness parameter of Mat\'ern kernel.}}
    \label{tab:summary}
\end{table}

\noindent
To hightlight our originality and contribution, a comparison of our results with the state-of-the-art complexity analysis is given in Tab.~\ref{tab:summary}. As far as we know, our work is the first to give a unified general lower bound in terms of metric entropy. Interestingly, we also notice that the commonly seen $\Theta(1/\epsilon^2)$ term in noisy setting disappear in noiseless setting, which matches our intuition that estimating a point with Gaussian noise typically takes $\Theta(1/\epsilon^2)$ sample complexity. Specifically, our contributions include:
\begin{itemize}
    \item We introduce a new set of analysis techniques and derive a \rev{\emph{general}} \rev{unified lower bound for} the deterministic problem complexity of \ego in terms of the \emph{metric entropy} of the function space ball in the corresponding reproducing kernel Hilbert space, \rev{providing a unified and intuitive understanding for the complexity of \ego}.    
    \item  Our \emph{general} lower bound allows us to leverage existing estimates of the covering number of the function space ball in the RKHS to derive kernel-specific lower bounds for the commonly used squared exponential kernel and Mat\'ern kernel with a large smoothness parameter $\nu$, \rev{without the commonly seen $1/\epsilon^2$ term for the noisy setting interestingly}. \rev{Furthermore, the lower bound for squared exponential kernel under noiseless setting is derived for the first time, to the best of our knowledge.} 
    \item We further show that these kernel-specific lower bounds nearly match the upper bounds attained by some non-adaptive search algorithms. Hence, our general lower bound is close to optimal for these specific kernels. 
\end{itemize}

\section{Related work}
There has been a large body of literature on analyzing the convergence properties of \ego. We first summarize the relevant literature area by area. We then highlight the position and the original contribution of our paper.

\bigskip
\noindent
\textbf{Algorithm-dependent Convergence Analysis.} One line of research analyzes the property of particular types of algorithms. For example, some papers~\cite{locatelli1997bayesian,frazier2011consistency} analyze the consistency of \ego{ algorithms}. \cite{vazquez2010convergence,wang2014theoretical} analyze the convergence property of the expected improvement algorithm. \cite{vakili2021optimal} proposes a maximum variance reduction algorithm that achieves optimal order simple regret for particular kernel functions. Under the assumption of H\"older continuity of the covariance function, lower and upper bounds are derived for the Bayesian setting in~\cite{grunewalder2010regret}. \rev{Among this set of literature, the works on information-theoretic upper bounds are more relevant to our metric-entropy lower bound.} \cite{srinivas2012information} derives an information-theoretic upper bound for the cumulative regret of the upper confidence bound algorithm.  \cite{russo2016information} gives an information-theoretic analysis of Thompson sampling. \rev{However, there is no existing work that provides a complementary information-theoretic lower bounds.} 

\bigskip
\noindent
\textbf{Kernel-specific Lower Bound Analysis.} As for lower bounds or complexity analysis, \cite{bull2011convergence} derives a lower bound of simple regret for Mat\'ern kernel in a noise-free setting. \cite{scarlett2017lower} provides lower bounds of both simple regret and cumulative regret for the squared exponential and Mat\'ern kernels. With the Mat\'ern kernel, a tight regret bound has been provided for Bayesian optimization in one dimension in~\cite{scarlett2018tight}. With heavy-tailed noise in the non-Bayesian setting, a cumulative regret lower bound has been provided for the Mat\'ern and squared exponential kernels in~\cite{ray2019bayesian}. More recently,~\cite{cai2021lower} provides lower bounds for both standard and robust Gaussian process bandit optimization. However, unlike the information-theoretic upper bound shown in~\cite{srinivas2012information}, the existing lower bound results are mostly~(if not all) restricted to specific kernel functions~(mostly squared exponential and Mat\'ern). The explicit connection between the optimization lower bound and the complexity of the RKHS has not been established so far in the existing literature. In this paper, we establish such a connection by constructing a lower bound in terms of \emph{metric entropy}.

\bigskip
\noindent
\textbf{Covering Number Estimate in RKHS.} Another area of research relevant to this paper is the estimate of covering number or metric entropy in function spaces. Some of the \rev{classical} results are used in this paper. \rev{In~\cite[Sec.~3.3]{edmunds1996function}, the covering number for the function space ball in a \rev{Besov} space is estimated.} A technique to derive a lower estimate of the covering number for a stationary kernel is developed in~\cite{zhou2003capacity}, and as an application, a lower bound of a function space ball's covering number for squared exponential kernel is derived.  

\bigskip
\noindent
\textbf{General Information-based Complexity Analysis.} Our focus is \ego{ in} this paper, due to its increasing popularity and lack of a unified and intuitive understanding for its complexity. Nevertheless, there has also been many classical works in the general area of information-based complexity analysis. For example, it is shown that the optimal convergence rates of global optimization is equivalent to those of approximation in the sup-norm~\cite{novak2008tractability}. However, approximation in the sup-norm itself is another hard problem with its complexity to be understood. There are also another set of results that try to connect the finite rank approximation, which is more general than sample based interpolation, with metric entropy~\cite{edmunds1996function,steinwart2017short,lorentz1966metric}. However, they can not be directly applied to our \ego{ problem}, due to the general finite rank approximation definitions that are inconsistent with our sample based \ego{ setting}.     
 
\bigskip
\noindent
\textbf{Minimax Rates for Kernel Regression.} In learning theory, there are well-established results on covering number bound of learning error. Many existing works \cite{cucker2002mathematical,raskutti2012minimax} derive covering number bounds for the generalization error of learning problem with RKHS or more general hypothetical spaces. However, in a typical learning setting, the sample points and corresponding observations are assumed to be identically and independently distributed, with observations corrupted by noise. However, the setting we considered in this paper is an essentially different global optimization problem. Specifically, our goal is to identify a solution with desired level of optimality and the sample point can be adaptively selected.  

\section{Problem Statement}
We consider \ego{ in} a non-Bayesian setting~\cite{srinivas2012information}. Specifically, we optimize a deterministic function $f$ from a reproducing kernel Hilbert space~(RKHS) $\mathcal{H}$ with input space $\mathbb{R}^d$, where $d$ is the dimension. $\mathcal{H}$ is equipped with the reproducing kernel $k(\cdot, \cdot):\mathbb{R}^d\times \mathbb{R}^d\rev{\to}\mathbb{R}$. Let $\mathcal{X}\subset\mathbb{R}^d$ be the feasible set of the optimization. In the following, we will use $[n]$ to denote the set $\{1, 2, \cdots, n\}$. We assume that 
\begin{assump}\label{assump:support_set}
$\mathcal{X}$ is compact.
\end{assump}
Assumption~\ref{assump:support_set} is reasonable because in many applications~(e.g., continuous hyperparameter tuning) of \ego, we are able to restrict the optimization into certain ranges \rev{based on domain knowledge}.
Regarding the black-box function $f\in\mathcal{H}$ that we aim to optimize, we assume that,
\begin{assump}
\label{assump:bounded_norm}
$\|{f}\|_\mathcal{H}\leq R$, where $R$ is a positive real number and $\|\cdot\|_\Hil$ is the norm induced by the inner product associated with $\Hil$.
\end{assump}
Assumption~\ref{assump:bounded_norm} requires that the function to be optimized is \rev{regular} in the sense that it has bounded norm in the RKHS. 
\begin{assump}
\label{assump:kernel_bound}
$k(x_1,x_2)\leq 1, \forall x_1,x_2\in \mathcal{X}$ and $k(x_1, x_2)$ is continuous on $\mathbb{R}^d\times \mathbb{R}^d$. 
\end{assump}
Assumption~\ref{assump:kernel_bound} is a common assumption for analyzing the convergence and complexity of \ego. It holds for a large class of commonly used kernel functions~(e.g., Mat\'ern kernel and squared exponential kernel) after normalization.  

\bigskip
\noindent
Our problem~\footnote{In the Gaussian process bandit literature, the maximizition formulation is usually adopted, while in the global optimization literature, the minimization formulation is usually adopted. Here, we adopt the latter.} is formulated as
\begin{equation}
\min_{x\in \mathcal{X}}\quad  f(x). 
\label{eqn:formulation}
\end{equation}
We know that 
\[
f(x_1)-f(x_2)=\langle f, k(x_1,\cdot)-k(x_2, \cdot)\rangle\leq\norm{f}_\Hil\norm{k(x_1,\cdot)-k(x_2,\cdot)}_\Hil.
\] 
Hence, it can be shown under Assumptions~\ref{assump:bounded_norm} and \ref{assump:kernel_bound}, that $f$ is continuous and thus \eqref{eqn:formulation} has an optimal solution on the compact set $\mathcal{X}$.   
As in standard \ego, we restrict ourselves to the zero-order oracle case. That is, our algorithm can only query the function value $f(x)$ but not higher-order information at a point $x$ in each step. Based on the function evaluations before the current step, the algorithm sequentially decides the next point to sample. In this paper, we only consider oracle query~(namely, function evaluation) complexity without considering the complexity of solving auxiliary optimization problems in typical \ego{ algorithms}.

\bigskip
\noindent
In this paper, we focus on the performance metric of \emph{simple regret} $\sr$. 
\begin{definition}[Simple regret]
After $t$ function evaluations, simple regret $\sr\coloneqq\min_{\tau\in[t]}f(x_\tau)-\min_{x\in\mathcal{X}}f(x)$, where $[t]\coloneqq\{1, 2,\cdots,t\}$. 
\end{definition}
Note that in some of the literature, simple regret is also defined as $f(\hat{x}_t)-\min_{x\in\mathcal{X}}f(x)$ where $\hat{x}_t$ is one additional point reported after $t$ steps. Since we can always pay one more function evaluation for the reported point, this definition difference will not impact our convergence or complexity analysis.

\section{Preliminary}
To analyze the problem complexity of \ego, we need a metric to measure the complexity of the RKHS. As an extreme example, if we choose a linear kernel, the underlying function to be optimized is a linear function. Hence, we can reconstruct it after a finite number of steps and compute the optimum without any error. The covering number is such a widely used metric to measure the complexity of an RKHS~\cite{zhou2002covering}. To facilitate our discussion, we introduce some concepts about the complexity of function sets.

\bigskip
\noindent
Given a normed vector space $(V,\norm{\cdot})$ and a subset $\gss\subset V$, for $\epsilon>0$, we make the following complexity related definitions~\cite{wu2017lecture}. 
\begin{definition}[$\epsilon$-covering]
$\{v_1,\cdots,v_N\}$ is an $\epsilon$-covering of $\gss$ if 
\[
\gss\subset\cup_{i\in[N]}B_{\norm{\cdot}}(v_i,\epsilon),
\]
where $B_{\norm{\cdot}}(v_i,\epsilon)$ is the ball in $V$ centered at $v_i$ with radius $\epsilon$ with respect to the norm $\norm{\cdot}$.
\end{definition}

\begin{definition}[$\epsilon$-packing] $\{v_1,\cdots,v_N\}\subset \gss$ is an $\epsilon$-packing of $\gss$ if 
\[
\min_{i\neq j}\norm{v_i-v_j}>\epsilon.
\]
\end{definition}
\begin{definition}[Covering number]
The covering number $\mathcal{N}(\gss, \epsilon, \norm{\cdot})$ is defined to be $\min\left\{n\,|\,\exists \epsilon\text{-covering } \{v_1,\cdots,v_n\} \text{ with cardinality } n\right\}$. 
\end{definition}

\begin{definition}[Packing number]
The packing number $\mathcal{M}(\gss, \epsilon, \norm{\cdot})$ is defined to be $\max\left\{n\,|\,\exists \epsilon\text{-packing } \{v_1,\cdots,v_n\} \text{ with cardinality } n\right\}$. 
\end{definition}

\begin{definition}[Metric entropy]
The metric entropy of $(\gss, \norm{\cdot})$ is defined to be $\log\mathcal{N}(\gss, \epsilon, \norm{\cdot})$, where $\mathcal{N}$ is the covering number. 
\end{definition}
It can be verified that,
\begin{proposition}[The.~\rom{4},~\cite{kolmogorov1959varepsilon}]
$\mathcal{N}(\gss,\epsilon, \norm{\cdot})\leq \mathcal{M}(\gss,\epsilon, \norm{\cdot})\leq \mathcal{N}(\gss,\frac{\epsilon}{2}, \norm{\cdot})$.
\end{proposition}
To facilitate the subsequent complexity analysis, we use $x_1, x_2,\cdots,x_t$ to denote the sequence of evaluated points up to step $t$. We now state the deterministic algorithm for solving the \ego{ problem}.

\begin{definition}[Deterministic algorithm]
A deterministic algorithm $\mathcal{A}$ for solving the optimization problem in~\eqref{eqn:formulation} is a sequence of mappings $(\pi_t)_{t=1}^\infty$, where $\pi_t:(\mathcal{X}\times\mathbb{R})^{t-1}\rev{\to}\mathcal{X},t\geq2$ and $\pi_1: \{\emptyset\}\rev{\to}\mathcal{X}$. When running the algorithm $\mathcal{A}$, the sample at step $t$ is $x_t=\pi_t((x_{\tau}, f(x_\tau))_{\tau=1}^{t-1}),t\geq2$ and $x_1=\pi_1(\emptyset)$. 
\end{definition}
\rev{Note that deterministic algorithms include most of the popular acquisition functions based \ego{ algorithms}~(e.g., lower/upper confidence bound~\cite{srinivas2012information} and expected improvement~\cite{jones1998efficient}).}

\bigskip
\noindent
We assume that the first sample point $x_1$ is deterministic, either given before running the algorithm or chosen by the algorithm. Now, if we suppose that $f$ is such that the algorithm observes a sequence of $0$'s for every function evaluation $f(x_\tau)$, it will generate a deterministic sample trajectory. We will see in our main result that this trajectory can be used to construct adversarial functions to derive the lower bound. We formally define it below. 
\begin{definition}[Zero sequence]
Given a deterministic algorithm $\mathcal{A}=(\pi_t)^\infty_{t=1}$. We set $x^0_1=\pi_1(\emptyset)$. Applying the recurrence relationship $x^0_t=\pi_t((x_\tau^0, 0)_{\tau=1}^{t-1})$, we get a deterministic sequence $x_1^0, x_2^0,\cdots, x_t^0,\cdots$, which only depends on the algorithm $\mathcal{A}$. We call this sequence the zero sequence of the algorithm $\mathcal{A}$.\label{def:zero_seq} 
\end{definition}

\section{Main Results}
Our strategy to derive the lower bound is decomposing the RKHS into two orthogonal subspaces with one of them expanding as more samples \rev{are} obtained. Then, we can project the function space ball into these two subspaces. We will show that as the number of sampled points grows, the covering number of the ball's projection into one subspace increases and the other decreases. We derive the lower bound on the number of optimization steps by bounding the increase/decrease rate. All the proofs of the lemmas and theorems are attached in the Appendix, except Lem.~\ref{lem:cover_ratio_bound} and Thm.~\ref{thm:lower_bound}. Before proceeding, we introduce some notations. 

\bigskip
\noindent
\textbf{Notations} For $f\in\Hil$, $f|_\mathcal{X}:\mathcal{X}\rev{\to}\mathbb{R}$ is defined as $f|_\mathcal{X}(x)=f(x),\forall x\in\mathcal{X}$. For $Q\subset\Hil$, we use $Q(\mathcal{X})$ to denote the set $\{f|_{\mathcal{X}}|f\in Q\}$, which is a subset of $C(\mathcal{X}, \norm{\cdot}_\infty)$, the continuous function space over $\mathcal{X}$, \rev{due to} Assumption~\ref{assump:kernel_bound}. $Q\rX$ is considered as a subset of $C(\mathcal{X}, \norm{\cdot}_\infty)$ in $\mathcal{N}(Q\rX, \epsilon, \norm{\cdot}_\infty)$ and $\mathcal{M}(Q\rX, \epsilon, \norm{\cdot}_\infty)$. 

\bigskip
\noindent
We first decompose the RKHS into two orthogonal subspaces.
\begin{definition}
$\Hil_t^{\|}:=\{\sum_{i\in[t]}\alpha_ik(x_i,\cdot)|\alpha_i\in\mathbb{R}\},$
$\Hil_t^{\perp}:=\{f\in\Hil|f(x_i)=0,\forall i\in[t]\}$.
\end{definition}
Notice that $\Hil_t^\|$ expands when we have more and more function evaluation data. In parallel, $\Hil_t^\perp$ shrinks.  
We then consider the intersection of the function space ball $S$ with $\Hil_t^\|$ and $\Hil_t^\perp$.
\begin{definition}
$S\coloneqq\{f|f\in\Hil, \norm{f}_\Hil\leq R\}, S_t^\|\coloneqq\Hil_t^\|\cap S, S_t^\perp\coloneqq\Hil_t^\perp\cap S$.
\end{definition}
With these definitions, we can show that any function in $S$ can be decomposed into two functions in $S_t^\|$ and $S_t^\perp$, respectively. 
\begin{restatable}{lemma}{lemfuncsplit}\label{lem:func_split}
$\forall f\in S$, there exists $m_t\in S_t^\|$, such that $f-{m}_t\in S_t^\perp$.
\end{restatable}
\begin{remark}
\label{rem:func_split}
 When the matrix $K=(k(x_i,x_j))_{i,j\in[t]}$ is invertible, we can check that ${m}_t(x)=f_X^TK^{-1}K_{Xx}$, where $f_X=[f(x_1), f(x_2),\cdots, f(x_t)]^T$ and $K_{Xx}=[k(x_1,x), k(x_2,x), \cdots, k(x_t, x)]^T$, satisfies $m_t\in S_t^\|$ and $f-{m}_t\in S_t^\perp$. The function $m_t(x)$ is exactly the posterior mean function in Gaussian process regression. 
\end{remark}
Intuitively, we can add some function from $S_t^\perp$ to $f$ without changing the historical evaluations at $x_1,\cdots,x_t$. If we have some way of lower bounding the complexity of $S_t^\perp$, we may be able to find a perturbing function from $S_t^\perp$ that leads to sub-optimality. We will try to lower bound the complexity of $S_t^\perp$ through Lem.~\ref{lem:decompose} and Lem.~\ref{lem:bound_para_covering}. 
\bigskip
\noindent
Since $S_t^\|$ and $S_t^\perp$ are orthogonal to each other in the RKHS, it is intuitive that the complexity of $S$ can be decomposed into the complexity of $S_t^\perp$ and $S_t^\|$. Formally, we have Lem.~\ref{lem:decompose}.
\begin{restatable}{lemma}{lemdecompose}\label{lem:decompose}
For any $\epsilon_t^\|>0, \epsilon_t^\perp>0$, we have 
\[
\mathcal{M}(S_t^\perp\rX, \epsilon_t^\perp, \norm{\cdot}_\infty)\geq\frac{\mathcal{N}(S\rX, \epsilon_t,\norm{\cdot}_\infty)}{\mathcal{N}(S_t^\|\rX, \epsilon_t^\|,\norm{\cdot}_\infty)},
\] 
where $\epsilon_t=\epsilon_t^\|+\epsilon_t^\perp$.
\end{restatable}
\noindent
Lem.~\ref{lem:decompose} is proved based on Lem.~\ref{lem:func_split}. With Lem.~\ref{lem:decompose}, we can lower bound $\mathcal{M}(S_t^\perp\rX, \epsilon_t^\perp, \norm{\cdot}_\infty)$ if we are able to upper bound $\mathcal{N}(S_t^\|\rX, \epsilon_t^\|,\norm{\cdot}_\infty)$.

\bigskip
\noindent
Since $S_t^\|$ is inside a finite dimensional space $\Hil_t^\|$, we can show that,
\begin{restatable}{lemma}{lemboundparacovering}
If $0<\epsilon<\frac{R}{4}$, we have $\log\mathcal{N}{(S_t^\|\rX, \epsilon,\norm{\cdot}_\infty)}\leq 2t\log\left(\frac{R}{\epsilon}\right)$. \label{lem:bound_para_covering} 
\end{restatable}
\noindent
We then give the following key lemma. 
\begin{lemma}
For $0<\epsilon<\epsilon_0$, if $t\leq\frac{\log\mathcal{N}(S\rX, 4\epsilon,\norm{\cdot}_\infty)}{4\log\left(\frac{R}{\epsilon}\right)}$, then for any sample sequence $x_1, \cdots, x_t$, we have, $$\frac{\mathcal{N}(S\rX, 4\epsilon,\norm{\cdot}_\infty)}{\mathcal{N}(S_t^\|\rX, \epsilon,\norm{\cdot}_\infty)}\geq2,$$ where $\epsilon_0=\min\left\{\sup\left\{\delta|\delta>0, \log\mathcal{N}(S\rX,4\delta, \norm{\cdot}_{\infty})>2\log2 \right\}, \frac{R}{4}\right\}$. 
\label{lem:cover_ratio_bound}
\end{lemma}
\Proof By assumption that $t\leq\frac{\log\mathcal{N}(S\rX, 4\epsilon,\norm{\cdot}_\infty)}{4\log\left(\frac{R}{\epsilon}\right)}$, we have $2t\log\left(\frac{R}{\epsilon}\right)\leq\frac{\log\mathcal{N}(S\rX, 4\epsilon,\norm{\cdot}_\infty)}{2}$. Meanwhile, by $\epsilon<\epsilon_0$ and the definition of $\epsilon_0$, $\frac{\log\mathcal{N}(S\rX, 4\epsilon,\norm{\cdot}_\infty)}{2}-\log2>0$. Combining with Lem.~\ref{lem:bound_para_covering}, we have,
\bee 
\begin{aligned}
&\log\mathcal{N}({S}_t^\|\rX, \epsilon,\norm{\cdot}_{\infty})\\
\leq&~2t\log\left(\frac{R}{\epsilon}\right)\\
\leq&~\frac{\log\mathcal{N}(S\rX, 4\epsilon,\norm{\cdot}_\infty)}{2}+\frac{\log\mathcal{N}(S\rX, 4\epsilon,\norm{\cdot}_\infty)}{2}-\log2\\
=&~\log\mathcal{N}(S\rX, 4\epsilon,\norm{\cdot}_\infty)-\log2.
\end{aligned}
\eee
So $\frac{\mathcal{N}(S\rX, 4\epsilon,\norm{\cdot}_\infty)}{\mathcal{N}(S_t^\|\rX, \epsilon,\norm{\cdot}_\infty)}\geq2$.
\qed

\bigskip
\noindent
We are now ready to give our main result in Thm.~\ref{thm:lower_bound}.
\begin{restatable}{theorem}{thmlowerbound}\label{thm:lower_bound}
If there exists a deterministic algorithm that achieves simple regret $\SR\leq\epsilon$ for any function $f\in S$ in $T$ function evaluations for our problem~\eqref{eqn:formulation}, it is necessary that,
\begin{equation}
    T=\Omega\left(\frac{\log\mathcal{N}(S\rX, 4\epsilon,\norm{\cdot}_\infty)}{\log(\frac{R}{\epsilon})}\right).
\end{equation}
\end{restatable}
\noindent
Before we prove Thm.~\ref{thm:lower_bound}, we give a sketch of the proof. 
For any deterministic algorithm and any number of optimization steps $t$, we consider the corresponding deterministic zero sequence $x^0_1, x^0_2, \cdots, x^0_t$ as defined in Def.~\ref{def:zero_seq}. We try to construct an adversarial function inside the corresponding $S_t^\perp$ with $0$ function value at the points $x^0_i,i\in[t]$ and low function values at some point that is not selected. The possible minimal value of such an adversarial function links to the covering number of the set $S_t^\perp\rX$, which can be lower bounded by combining Lem.~\ref{lem:decompose} and Lem.~\ref{lem:bound_para_covering}. 

\begin{proof}[\textbf{Proof of Thm.~\ref{thm:lower_bound}}]
Given an deterministic algorithm $\mathcal{A}=(\pi_t)_{t=1}^{+\infty}$, if it always gets the evaluations $0$, then the sample trajectory satisfies,
$$
x_t^0 = \pi_t\left((x_\tau^0,0)_{\tau=1}^{t-1}\right), t\geq2,
$$
which is exactly the zero sequence of the algorithm. Note that the zero sequence $x_t^0$ only depends on the deterministic algorithm $\mathcal{A}$. Once we fix the algorithm, the zero sequence is fixed.   

\bigskip
\noindent
We want to check the feasibility of the problem~\eqref{prob:bad_func},
\begin{equation}
 \begin{aligned}
\underset{s \in \mathcal{S}, x\in \mathcal{X}}{ \min } &~1 \\
\text { s.t. } & s\left(x_{n}^0\right)=0,~\forall n=1, \ldots, t,\\
&s(x)<-\epsilon. 
\end{aligned}  
\label{prob:bad_func}
\end{equation}
Any feasible solution of~\eqref{prob:bad_func} has some `adversarial' property against the algorithm $\mathcal{A}$. In fact, suppose that $(\bar{s},\bar{x})$ is a feasible solution for problem~\eqref{prob:bad_func}, when we run the algorithm $\mathcal{A}$ over $\bar{s}$, the sample sequence up to step $t$ is exactly the zero sequence truncated at step $t$ and $r_{(t)}=\min_{\tau\in[t]}\bar{s}(x_\tau^0)-\min_{x\in\mathcal{X}}\bar{s}(x)>\epsilon$. Now the question is under what condition, the problem~\eqref{prob:bad_func} is feasible. 
Since we are analyzing the asymptotic rate, we restrict to the case $\epsilon<\epsilon_0$, where $\epsilon_0$ is given in Lem.~\ref{lem:cover_ratio_bound}. 

\bigskip
\noindent
By Lem.~\ref{lem:cover_ratio_bound} and Lem.~\ref{lem:decompose}, if $t\leq\frac{\log\mathcal{N}(S\rX, 4\epsilon,\norm{\cdot}_\infty)}{4\log\left(\frac{R}{\epsilon}\right)}$, for the sample sequence $x_1^0, \cdots, x_t^0$ corresponding to any given algorithm, we have, 
$$
\mathcal{M}(S_t^\perp\rX, 3\epsilon, \norm{\cdot}_\infty)\geq\frac{\mathcal{N}(S\rX, 4\epsilon,\norm{\cdot}_\infty)}{\mathcal{N}(S_t^\|\rX, \epsilon,\norm{\cdot}_\infty)}\geq 2.
$$
Therefore, there exists functions $f_1, f_2\in S_t^\perp$, such that, $\norm{f_1|_\mathcal{X}-f_2|_\mathcal{X}}_\infty\geq 3\epsilon$. So $\norm{f_1|_\mathcal{X}}_\infty+\norm{f_2|_\mathcal{X}}_\infty\geq3\epsilon$ and at least one of $f_1$ and $f_2$ has \rev{$L_\infty$} norm over the set $\mathcal{X}$ at least $\frac{3\epsilon}{2}$. Without loss of generality, we assume $\norm{f_1|_\mathcal{X}}_\infty\geq\frac{3\epsilon}{2}$. Since for $\forall g\in S_t^\perp$, $-g\in S_t^\perp$, there exists $\hat{f}\in S_t^\perp$~(either $f_1$ or $-f_1$), such that, $$\inf_{x\in\mathcal{X}}\hat{f}(x)\leq -\frac{3\epsilon}{2}.$$ When applying the given algorithm to $\hat{f}$, if $t\leq\frac{\log\mathcal{N}(S\rX, 4\epsilon,\norm{\cdot}_\infty)}{4\log\left(\frac{R}{\epsilon}\right)}$, the suboptimality gap or the simple regret $r_{(t)}$ is at least $\frac{3}{2}\epsilon$. Therefore, to reduce the simple regret $r_{(T)}\leq\epsilon$ for all the functions in $S$ within $T$ steps, it is necessary that, 
\begin{equation}
    T=\Omega\left(\frac{\log\mathcal{N}(S\rX, 4\epsilon,\norm{\cdot}_\infty)}{\log(\frac{R}{\epsilon})}\right).
\end{equation}
\qed
\end{proof}
To verify the effectiveness of Thm.~\ref{thm:lower_bound}, we apply it to a simple case in Ex.~\ref{ex:simple_ex_verify}. 
\begin{example}\label{ex:simple_ex_verify}
For the quadratic kernel $k(x,y)=(x^Ty)^2$, the corresponding RKHS is finite dimensional and is given as~\cite{mairal2018machine},
\begin{equation}
\Hil=\left\{f_A(x)=x^TAx|A\in\mathcal{S}^{d\times d}\right\},    
\end{equation}
where $\mathcal{S}^{d\times d}$ is the set of symmetric matrices of size $d\times d$. We know that,
\begin{equation}
    \langle f_{A_1}, f_{A_2}\rangle_\Hil = \langle A_1, A_2\rangle_\mathrm{F}, 
\end{equation}
where $\langle\cdot,\cdot\rangle_\mathrm{F}$ is the Frobenius inner product. Since $\mathcal{S}^{d\times d}$ can be embedded into $\mathbb{R}^{\frac{d\times(d+1)}{2}}$ and the metric entropy for compact set in Euclidean space is ${\Theta}\left(\log\frac{1}{\epsilon}\right)$~\cite{wu2017lecture}, the lower bound in Thm.~\ref{thm:lower_bound} reduces to a constant. By applying a grid search algorithm for the quadratic kernel, we can identify the ground truth function after a finite number of steps and determine the optimal solution without any error. Therefore, the lower bound is tight \rev{in $\epsilon$} in the case of the quadratic kernel.
\end{example}

\subsection{Comparison with upper bounds for commonly used kernels}
Ex.~\ref{ex:simple_ex_verify} demonstrates the validity of Thm.~\ref{thm:lower_bound} for simple quadratic kernel functions. In this section, we will derive kernel-specific lower bounds for the squared exponential kernel and the Mat\'ern kernels by using existing estimates of the covering numbers for their RKHS's and Thm.~\ref{thm:lower_bound}. We compare our lower bounds with derived/existing upper bounds and show that they nearly match.     

\subsubsection{Squared Exponential kernel}
One widely used kernel in \ego{ is} the squared exponential~(SE) kernel given by
\begin{equation}
   k(x,y) = \exp{\left\{-\frac{\norm{x-y}^2}{\sigma^2}\right\}}. \label{eq:se} 
\end{equation}
In this case, we restrict to $\mathcal{X}=[0,1]^d$. By applying Thm.~\ref{thm:lower_bound}, we have,
\begin{restatable}{theorem}{thmselbub}
With $\mathcal{X}=[0,1]^d$ and using the squared exponential kernel, if there exists a deterministic algorithm that achieves simple regret $\SR\leq\epsilon$ for any function $f\in S$ in $T$ function evaluations for our problem~\eqref{eqn:formulation}, it is necessary that,
\begin{equation}
    T=\Omega\left(\left(\log\frac{R}{\epsilon}\right)^{d/2-1}\right).
\end{equation}
Furthermore, there exists a deterministic algorithm and $T$ satisfying, 
$$
T=\mathcal{O}\left(\left(\log\frac{R}{\epsilon}\right)^{d}\right),
$$
such that the algorithm achieves $\SR\leq\epsilon$ in $T$ function evaluations for any $f\in S$.
\label{thm:se_lb_ub}
\end{restatable}
\noindent
The upper bound part is obtained through sampling non-adaptively to reduce the posterior variance to a uniform low level in $\mathcal{X}$. In this theorem, we focus on the asymptotic analysis of \ego{ and} hide the coefficients that may depend on the dimension. We notice that the upper bound and lower bound are both polynomial in $\log\frac{1}{\epsilon}$. The orders of the polynomials nearly match, up to a replacement of $d/2$ by $d$ and a logarithmic term $\log\frac{R}{\epsilon}$.    

\subsubsection{Mat\'ern kernel}\label{subsubsec:mat_res}
In this section, we consider the Mat\'ern kernel,
\begin{equation}
 k(x,y)=C_{\nu}(\norm{x-y})=\sigma^{2} \frac{2^{1-\nu}}{\Gamma(\nu)}\left(\sqrt{2 \nu} \frac{\norm{x-y}}{\rho}\right)^{\nu} K_{\nu}\left(\sqrt{2 \nu} \frac{\norm{x-y}}{\rho}\right),
\end{equation}
where $\rho$ and $\nu$ are positive parameters of the covariance, $\Gamma$ is the gamma function, and $K_{\nu}$ is the modified Bessel function of the second kind.


\begin{restatable}{theorem}{thmmatlbub}
With $\mathcal{X}=[0, 1]^d$ and the Mat\'ern kernel, if there exists a deterministic algorithm that achieves simple regret $\SR\leq\epsilon$ for any function $f\in S$ in $T$ function evaluations for our problem~\eqref{eqn:formulation}, it is necessary that,
\begin{equation}
 T=\Omega\left(\left(\frac{R}{\epsilon}\right)^{\frac{d}{\nu+d/2}}\left(\log\frac{R}{\epsilon}\right)^{-1}\right). 
\end{equation}
Furthermore, there exists a deterministic algorithm and $T$ satisfying, 
\begin{equation}
 T=\mathcal{O}\left(\left(\frac{R}{\epsilon}\right)^{\frac{d}{\nu}}\right),\label{eq:mat_upper_bound} 
\end{equation}
such that the algorithm achieves $\SR\leq\epsilon$ in $T$ function evaluations for any $f\in S$.
\label{thm:mat_lb_ub}
\end{restatable}


\begin{remark}
The upper bound part of Thm.~\ref{thm:mat_lb_ub} is proved by Thm.~1 of \cite{bull2011convergence}. We also notice that \cite{bull2011convergence} provides a lower bound of the same order as the upper bound in Eq.~\eqref{eq:mat_upper_bound}, which means that the upper bound order is also the \rev{optimal} lower bound order. 
\end{remark}
\begin{remark}\label{remark:mat_low}
When ${\nu}\geq\frac{1}{2}d$, our lower bound can further imply the lower bound of $\Omega\left(\left(\frac{R}{\epsilon}\right)^{\frac{d}{2\nu}}\left(\log\frac{R}{\epsilon}\right)^{-1}\right)$, which nearly matches the upper bound, up to a replacement of $d/2$ by $d$ and a logarithmic term $\log\frac{R}{\epsilon}$. However, when $\frac{\nu}{d}$ is small, there is still a significant gap between the lower bound implied by our general lower bound and the optimal lower bound.
\end{remark}
\subsection{Discussions and limitations}\label{subsec:discuss_limit}
\subsubsection{Comparison with the noisy case}
In~\cite{scarlett2017lower}, lower bounds on simple regret of \ego{ with} Mat\'ern and squared exponential kernels in a noisy setting were provided. Specifically, a lower bound of $\Omega\left(\frac{1}{\epsilon^{2}}\left(\log \frac{1}{\epsilon}\right)^{d / 2}\right)$ for the squared exponential kernel and $\Omega\left(\left(\frac{1}{\epsilon}\right)^{2+d / \nu}\right)$ for the Mat\'ern kernel are provided in~\cite{scarlett2017lower}. We notice that there is one common $\left(\frac{1}{\epsilon}\right)^2$ term for both kernels, in addition to the lower bound in the noiseless case. In particular, for the squared exponential kernel, this term dominates the lower bound. \rev{Interestingly, this $\Theta(1/\epsilon^2)$ term in noisy setting disappear in noiseless setting, which matches our intuition that estimating a point with Gaussian noise typically takes $\Theta(1/\epsilon^2)$ sample complexity.} We thus conjecture that, $\Omega\left(\frac{1}{\epsilon^2}\frac{\log\mathcal{N}(S\rX, 4\epsilon,\|\cdot\|_\infty)}{\log(\frac{1}{\epsilon})}\right)$ is a general lower bound for the noisy case, though we do not currently have a proof. 


\subsubsection{Connection with cumulative regret}
In some applications, the black-box optimization problem is solved in a finite number of steps and then a fixed solution is reported and applied. In such scenarios, simple regret or suboptimality of the solution finally reported is of interest. In contrast, in the bandit settings~\cite{agrawal1995continuum}, we are typically interested in the cumulative cost over a horizon, where cumulative regret, \rev{defined} as the accumulation of the suboptimality gap of the samples, will be more relevant. In our main results, we did not directly analyze the cumulative regret. But if we assume that the total number of function evaluations $T$ is revealed beforehand, we can derive bounds on cumulative regret using the upper bound parts of Thm.~\ref{thm:se_lb_ub} and Thm.~\ref{thm:mat_lb_ub}. For the squared exponential kernel, by allocating the first $\Theta\left((\log T)^d\right)$ number of evaluations to reduce the suboptimality gap to the level of $\mathcal{O}\left(\frac{1}{T}\right)$ and continue to evaluate the solution with at most $\mathcal{O}\left(\frac{1}{T}\right)$ suboptimality gap later, the total cumulative regret is upper bounded by $\mathcal{O}((\log T)^d)$. Similarly, for the Mat\'ern kernel, by allocating the first $\Theta\left(T^{\frac{d}{d+\nu}}\right)$ number of evaluations to reduce the suboptimality gap to the level of $\mathcal{O}\left(\left(\frac{1}{T}\right)^{\frac{\nu}{d+\nu}}\right)$ and continue to evaluate the solution with at most $\mathcal{O}\left(\left(\frac{1}{T}\right)^{\frac{\nu}{d+\nu}}\right)$ suboptimality gap later, the total cumulative regret is upper bounded by $\mathcal{O}(T^{\frac{d}{d+\nu}})$. \rev{Standard ``doubling trick'' can be used to design an algorithm with cumulative regret bound of the same order without knowing $T$ beforehand}. 

\section{Experiments} \label{sec:exp}
In this section, we will first give a demonstration of adversarial functions, on which two common algorithms, the lower confidence bound~(LCB) and the expected improvement~(EI), perform poorly and achieve the optimization lower bound. Then we run the two algorithms on a set of randomly sampled functions and compare the average performance and the adversarial performance in terms of simple regret. The algorithms are implemented based on~\textsf{GPy}~\cite{gpy2014} and \textsf{CasADi}~\cite{andersson2019casadi}. All the auxiliary optimization problems in the algorithms are solved using the solver \textsf{IPOPT}~\cite{wachter2006implementation} with multiple different starting points. Our experiments take about 15 hours on a device with \texttt{AMD Ryzen Threadripper 3990X 64-Core Processor} and 251 GB RAM. 
 
\subsection{Demonstration of adversarial functions}
 \begin{figure}[htbp]
    \centering
    \includegraphics[width=\linewidth]{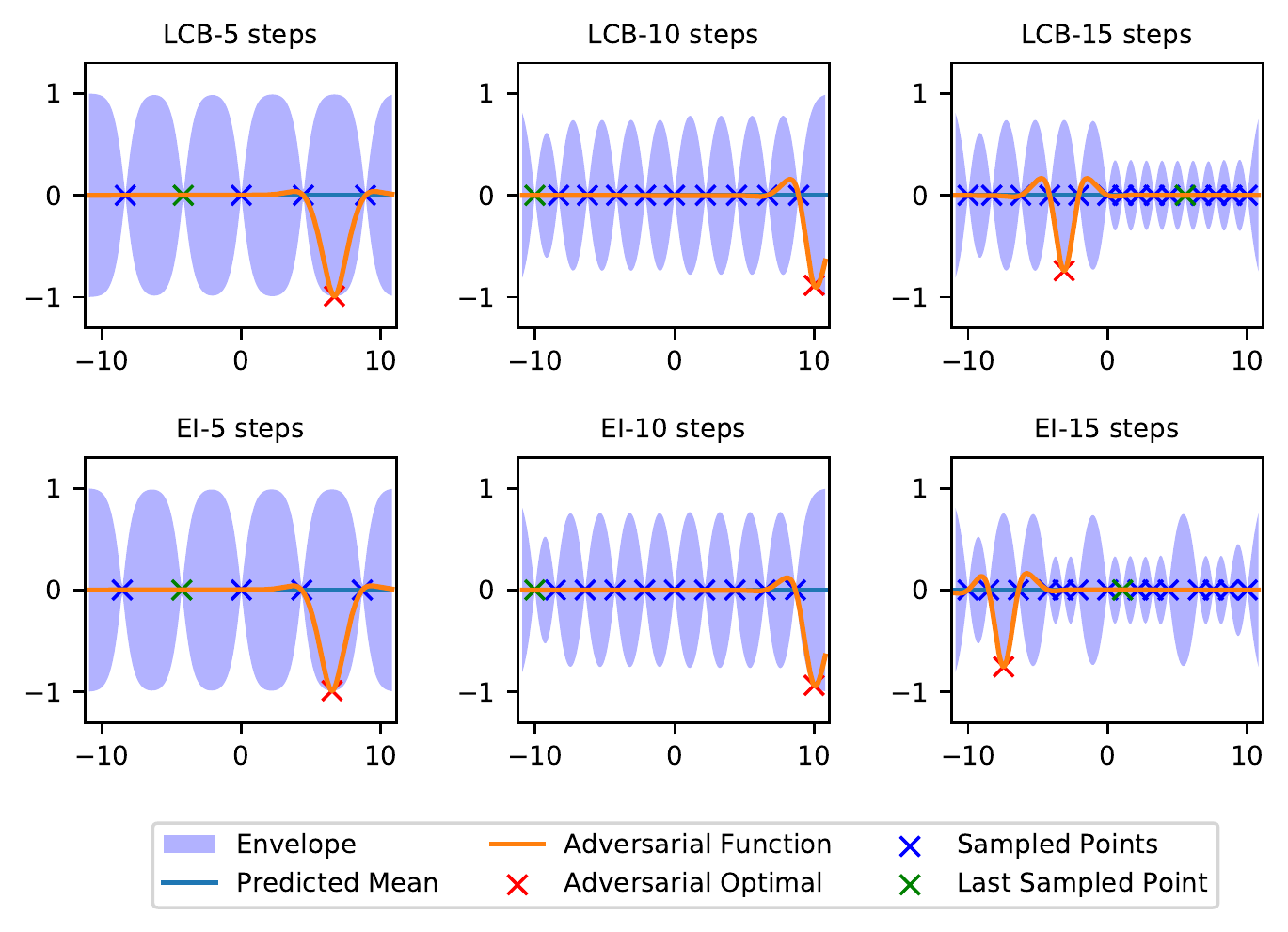}
    \caption{Demonstrations of adversarial functions in dimension one.}
    \label{fig:adv_func_1d}
\end{figure}
\noindent
In our proof of Thm.~\ref{thm:lower_bound}, we use a particular set of adversarial functions, which reveal value $0$ to the algorithm and have low values somewhere else. In this section, we demonstrate such adversarial functions for two popular algorithms, expected improvement and lower confidence bound. 

\bigskip
\noindent
We use the Mat\'ern kernel in one dimension with $\nu=\frac{5}{2}, \rho=1, \sigma^2=1$. We set the compact set to $\mathcal{X}=[-10, 10]$ and assume that the RKHS norm upper bound is $R=1$. We apply both lower confidence bound algorithm with the constant weight $1$ for the posterior standard deviation and the expected improvement algorithm. We manually assign $x_1=0$ as the first sampled point and derive the adversarial function by solving Prob.~\eqref{opt:get_adv_func}.   
\begin{equation}
    \begin{aligned}
\min_{x\in\mathcal{X}}\min_{s \in \mathcal{H}}& ~s(x) \\
\text { s.t. } & s\left(x_{n}^0\right)=0, \forall n=1, \ldots, t \\
&\norm{s}_\Hil\leq R \\ 
\end{aligned}
\label{opt:get_adv_func}
\end{equation}
Thanks to the \rev{optimal recovery property~\cite[Thm 13.2]{wendland2004scattered}}, the optimal value for the inner problem \rev{of \eqref{opt:get_adv_func}} can be analytically derived as $$-R\sqrt{k(x,x)-k(x, X)^TK^{-1}k(X,x)}.$$
Fig.~\ref{fig:adv_func_1d} demonstrates the adversarial functions inside the corresponding RKHS with bounded norm of $1$, which have value $0$ at all the sampled points but have low global optimal value somewhere else. We notice that the envelope formed by the functions inside the ball with consistent evaluation data shrinks as more and more data becomes available. Intuitively, any algorithm needs to sample sufficiently densely globally in the adversarial case in order to find a close-to-optimal solution.  

\subsection{Average vs. adversarial performance}
\begin{figure}[htbp]
    \centering
    \includegraphics[width=\linewidth]{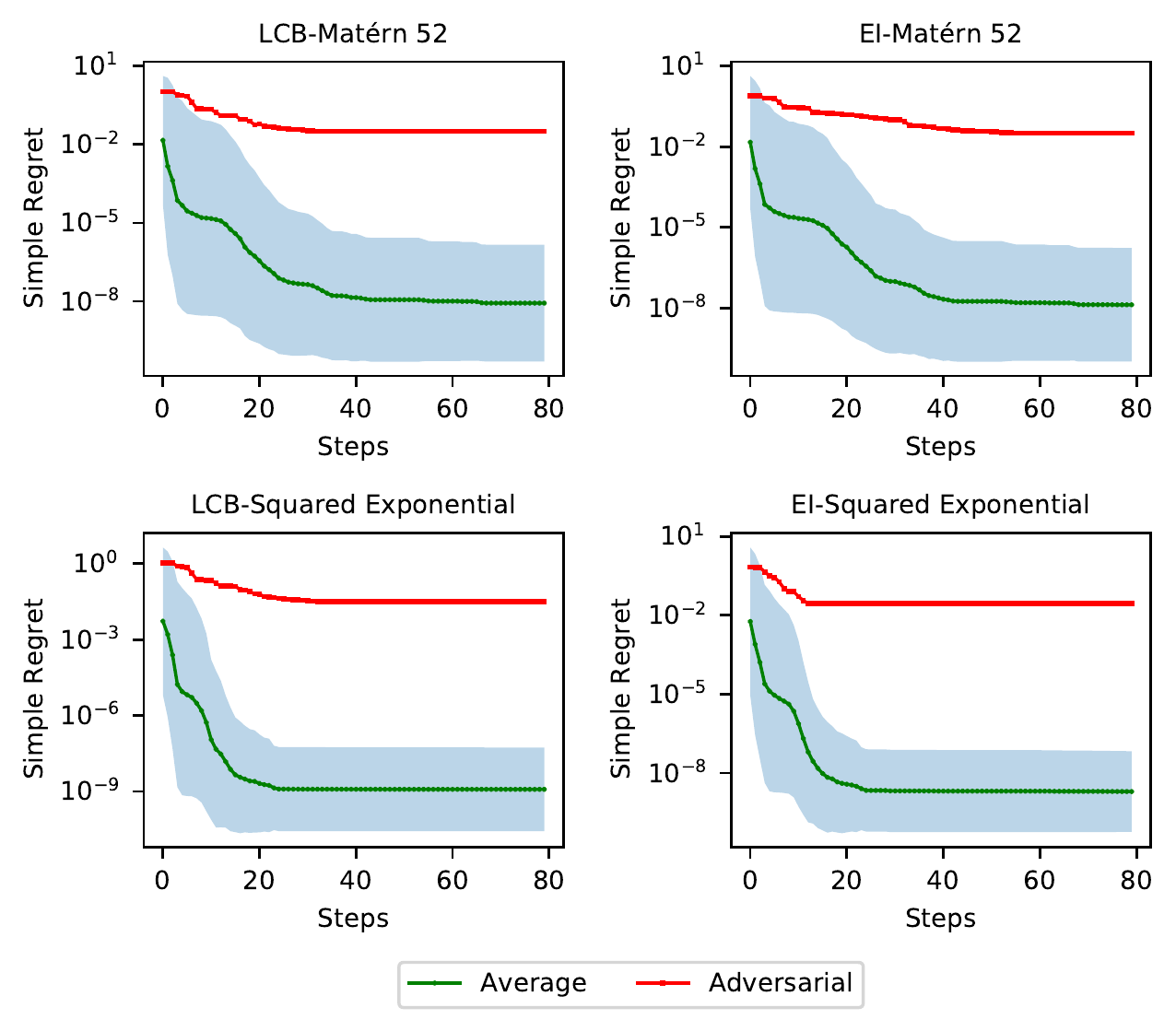}
    \caption{Comparison of average performance~($\pm$standard deviation shown as shaded area, over $100$ instances) and adversarial performance. Adversarial simple regret is defined as opposite of the optimal value of Prob.~\eqref{opt:get_adv_func}, namely the simple regret of the adversarial function at different optimization steps.}
    \label{fig:ave_adv_comp}
\end{figure}

The proofs of Thm.~\ref{thm:se_lb_ub} and Thm.~\ref{thm:mat_lb_ub} indicate that a non-adaptive sampling algorithm can achieve a close-to-optimal worst-case convergence rate. However, in practice, adaptive algorithms~(e.g., lower confidence bound and expected improvement) are usually adopted and perform better. There could potentially be a gap between average-case convergence and worst-case convergence. To perform such a comparison, we randomly sample a set of functions from the RKHS to run the algorithms over. Specifically, we first uniformly sample a finite set of knots $X\subset\mathcal{X}$ and then sample the function values $f_X$ on the knots from the marginal distribution of the Gaussian process, which is a finite-dimensional Gaussian distribution. We then construct the minimal norm interpolant of the knots as the sampled function. To be consistent with the bounded norm assumption, we reject the functions with a norm value larger than $R$.    

\bigskip
\noindent
We use simple regret, which is defined to be $\min_{\tau\in[t]}f(x_\tau)-\min_{x\in\mathcal{X}}f(x)$, to measure the performance of different algorithms. We set $\mathcal{X}=[0, 1]^3\subset\mathbb{R}^3$ and set the length scales and variances of both the Mat\'ern kernel function~(with $\nu=\frac{5}{2}$) and the squared exponential kernel. Fig.~\ref{fig:ave_adv_comp} shows the comparison of average simple regret and adversarial simple regret. We observe that the average performance is much better than the performance on adversarial functions in terms of simple regret. Intuitively, adversarial functions are only a subset of \emph{needle-in-haystack} functions, with most region flat and somewhere very small, when $t$ becomes large. For those adversarial functions such as shown in Fig.~\ref{fig:adv_func_1d}, it can be difficult for the \ego{ algorithms} to ``see'' the trend of the function. For common functions inside the function space ball, however, the algorithms are still able to detect the trend of the function value and find a near-to-optimal solution quickly. 

\section{Conclusions}
In this paper, we provide a general lower bound on the worst-case suboptimality or simple regret for noiseless \ego{ in} a non-Bayesian setting in terms of the metric entropy of the corresponding reproducing kernel Hilbert space~(RKHS). We apply the general lower bounds to commonly used specific kernel functions, including the squared exponential kernel and the Mat\'ern kernel. We further derive upper bounds and compare them to the lower bounds and find that they nearly match, except for the case for the Mat\'ern kernel when $\frac{\nu}{d}$ is small. Two interesting future research directions are deriving an upper bound on the worst-case convergence rate in terms of metric entropy and characterizing the average-case convergence rate. We also conjecture that introducing randomness into the existing algorithms can improve the worst-case performance. A expected analysis challenge is that our current approach is sensitive to randomness.

\begin{acknowledgements}
    This work was supported by the Swiss National Science Foundation under the NCCR Automation project, grant agreement 51NF40\_180545 and the RISK project (Risk Aware Data-Driven Demand Response, grant number 200021\_175627).
\end{acknowledgements}

\appendix  
\section*{Appendix}
In the appendices, we give proofs of all the lemmas and theorems in the paper. We also include information on computational resources and time. \section{Proof of Lemma~\ref{lem:func_split}}
\lemfuncsplit*
\label{sec:pf_lem_func_split}
\Proof
Consider the optimization problem below, 
\begin{equation}
    \begin{aligned}
\underset{s \in \mathcal{H}}{ \min } &~\|s\|_{\mathcal{H}}^{2} \\
\text { s.t. } & s\left(x_{n}\right)=f(x_{n}) \\
& \forall n=1, \ldots, t
\end{aligned}
\label{opt:get_m}
\end{equation}
Based on the representer theorem~\cite[Theorem 1.3.1]{wahba1990spline}, the optimal solution of \eqref{opt:get_m} has the form $\sum_{i=1}^t\alpha_ik(x_i, \cdot)$. By using the constraint $s(x_n)=f(x_n)$, we can derive $K\alpha=f_X$, where $f_X=[f(x_1), f(x_2), \cdots, f(x_n)]^T$ and $K=(k(x_i,x_j))_{i\in[n], j\in[n]}$. With this restriction, we transform the problem in~\eqref{opt:get_m} to the problem in~\eqref{opt:Ka_get_m}.
\begin{equation}
    \begin{aligned}
\underset{\alpha \in \mathbb{R}^t}{ \min } &~\alpha^TK\alpha \\
\text { s.t. } & K\alpha=f_X\\
\end{aligned}
\label{opt:Ka_get_m}
\end{equation}
We take $\alpha^*$ as the solution to the problem in~\eqref{opt:Ka_get_m}, whose feasibility is guaranteed by representer theorem~\cite{wahba1990spline} and the non-emptiness of the feasible set~($f$ is feasible for \eqref{opt:get_m}). Therefore, $m_t(x)=(\alpha^*)^TK_{Xx}$ is the optimal solution to \eqref{opt:get_m}. Since $f$ is a feasible solution for the problem~\eqref{opt:get_m}, $\norm{m_t}_\Hil\leq\norm{f}_\Hil\leq R$. In addition, $(f-m_t)(x_i)=f(x_i)-m_t(x_i)=0, \forall i\in[t]$. And $$
\begin{aligned}
\norm{f-m_t}_\Hil^2&=\norm{f}_\Hil^2+\norm{m_t}_\Hil^2-2\langle f, m_t\rangle\\
&=\norm{f}_\Hil^2-\norm{m_t}_\Hil^2\\
&\leq R^2.
\end{aligned}
$$
So $m_t\in S_t^\|$ and $f-m_t\in S_t^\perp$.
\qed

\section{Proof of Lemma~\ref{lem:decompose}}
\label{sec:pf_lem_decompose}
\lemdecompose*
\Proof
Let $(p_1,p_2,\cdots,p_m)$ be an $\epsilon_t^\|$-covering of $S_t^\|\rX$ and $(q_1,q_2,\cdots,q_n)$ an $\epsilon_t^\perp$-covering of $S_t^\perp\rX$. Then $\forall f\in S$, by Lem.~\ref{lem:func_split}, $f=m_t+(f-m_t)$, where $m_t\in S_t^\|$ and $f-m_t\in S_t^\perp$. By the definition of covering, $\exists p_i$, such that $\norm{m_t|_{\mathcal{X}}-p_i}_\infty\leq\epsilon_t^\|$ and $\exists q_j$, such that $\norm{(f-m_t)|_{\mathcal{X}}-q_j}_\infty\leq\epsilon_t^\perp$. So 
\[
\begin{aligned}
\norm{f|_\mathcal{X}-(p_i+q_j)}_\infty\leq&\norm{m_t|_\mathcal{X}-p_i}_\infty+\norm{(f-m_t)|_\mathcal{X}-q_j}_\infty\\
\leq&\epsilon_t^\|+\epsilon_t^\perp=\epsilon_t.
\end{aligned}
\] 
So the set $\{p_i+q_j|i\in[m],j\in[n]\}$ is an $\epsilon_t$-covering of $S\rX$ and we have the cardinality \[
|\{p_i+q_j|i\in[m],j\in[n]\}|=\mathcal{N}(S_t^\perp\rX, \epsilon_t^\perp, \norm{\cdot}_\infty)\mathcal{N}(S_t^\|\rX, \epsilon_t^\|, \norm{\cdot}_\infty)\geq\mathcal{N}(S\rX, \epsilon_t,\norm{\cdot}_\infty).
\] So $\mathcal{M}(S_t^\perp\rX, \epsilon_t^\perp, \norm{\cdot}_\infty)\geq\mathcal{N}(S_t^\perp\rX, \epsilon_t^\perp, \norm{\cdot}_\infty) \geq\frac{\mathcal{N}(S\rX, \epsilon_t,\norm{\cdot}_\infty)}{\mathcal{N}(S_t^\|\rX, \epsilon_t^\|,\norm{\cdot}_\infty)}$.
\qed

\section{Proof of Lemma~\ref{lem:bound_para_covering}}
\label{sec:pf_lem_bound_para_covering}
\lemboundparacovering*
\Proof
We first introduce the set,
$$
\mathcal{E}_t=\{\alpha\in\mathbb{R}^t|\alpha^T K_t\alpha\leq R^2\},
$$
where $K_t=(k(x_i,x_j))_{i,j\in[t]}$. Without loss of generality, we assume that $K_t$ has full rank in the following analysis. Notice that if this condition does not hold, we only need to restrict to the subspace spanned by the eigenvectors of $K_t$ with strictly positive eigenvalues and consider the intersection of $\mathcal{E}_t$ with the subspace. Since the restriction only reduces the essential dimension, the upper bound still holds. We introduce the norm $\norm{\alpha}_{K_t}=\sqrt{\alpha^TK_t\alpha}$.
We then have,
$\forall f(x)=\alpha^Tk(X,x)\in S_t^\|\rX, g(x)=\beta^Tk(X,x)\in S_t^\|\rX$, we have  
\bee
\begin{aligned}
&\norm{f-g}_\infty\\
=&~\sup_{x\in \mathcal{X}}|(\alpha-\beta)^Tk(X,x)|\\
=&~\sup_{x\in \mathcal{X}}|\langle\sum_{i\in[t]}(\alpha_i-\beta_i)k(x_i,\cdot), k(x,\cdot)\rangle|\\
\leq&~\sup_{x\in \mathcal{X}}\norm{\sum_{i\in[t]}(\alpha_i-\beta_i)k(x_i,\cdot)}_\Hil\norm{k(x,\cdot)}_\Hil\\
\leq&~\sup_{x\in \mathcal{X}}\norm{\alpha-\beta}_{K_t}\sqrt{k(x,x)}\\
\leq&\norm{\alpha-\beta}_{K_t}
\end{aligned}
\eee
Therefore, we have $\mathcal{N}(S^\|_t\rX, \epsilon,\norm{\cdot}_\infty)\leq\mathcal{N}(\mathcal{E}_t, \epsilon,\norm{\cdot}_{K_t})$.
We further have,
\bee 
\begin{aligned}
&\mathcal{N}(\mathcal{E}_t, \epsilon,\norm{\cdot}_{K_t})\\
\leq&~\mathcal{M}(\mathcal{E}_t, \epsilon,\norm{\cdot}_{K_t})\\
\leq&~\frac{\mathrm{Vol}\left(B_{\norm{\cdot}_{K_t}}\left(0, R+\frac{\epsilon}{2}\right)\right)}{\mathrm{Vol}\left(B_{\norm{\cdot}_{K_t}}\left(0, \frac{\epsilon}{2}\right)\right)}\\
=&~\left(\frac{2R}{\epsilon}+1\right)^t\\
\leq&~\left(\frac{R^2}{2\epsilon^2}+\frac{R^2}{2\epsilon^2}\right)^t\\
=&~\left(\frac{R}{\epsilon}\right)^{2t},
\end{aligned}
\label{eq:bound_inq}
\eee
The second inequality in~\eqref{eq:bound_inq} follows by that if $\alpha_1, \alpha_2,\cdots,\alpha_M$ is an $\epsilon$-packing of the set $\mathcal{E}_t$, then $\cup_{i\in[M]}B_{\norm{\cdot}_{K_t}}\left(\alpha_i, \frac{\epsilon}{2}\right)\subset B_{\norm{\cdot}_{K_t}}\left(0, R+\frac{\epsilon}{2}\right)$ and $B_{\norm{\cdot}_{K_t}}\left(\alpha_i, \frac{\epsilon}{2}\right)\cap B_{\norm{\cdot}_{K_t}}\left(\alpha_j, \frac{\epsilon}{2}\right)=\emptyset,\forall i\neq j$ by the definition of packing. The third inequality in~\eqref{eq:bound_inq} follows by the assumption of $0<\epsilon<\frac{R}{4}$. So $\log\mathcal{N}(\mathcal{E}_t, \epsilon,\norm{\cdot}_{K_t})\leq 2t\log\left(\frac{R}{\epsilon}\right)$. Therefore, $\log\mathcal{N}(S^\|_t\rX, \epsilon,\norm{\cdot}_\infty)\leq 2t\log\left(\frac{R}{\epsilon}\right)$. 
\qed


\section{Proof of Theorem~\ref{thm:se_lb_ub}}
\thmselbub*
\Proof
By~\cite[Example~1]{zhou2003capacity}, the covering number satisfies, 
\begin{equation}
   \log\mathcal{N}(S\rX, 4\epsilon, \norm{\cdot}_\infty)=\Omega\left(\log\left(\frac{R}{\epsilon}\right)^{\frac{d}{2}}\right). 
\end{equation}
Therefore, Thm.~\ref{thm:lower_bound} implies that,
\begin{equation}
   T = \Omega\left(\left(\log\frac{R}{\epsilon}\right)^{\frac{d}{2}-1}\right). 
\end{equation}
We now focus on proving the upper bound part. To facilitate the following proof, we define,
\begin{equation}
\begin{aligned}
\ubar{f}_t(x) &= m_t(x)-\sigma_t(x)\sqrt{R^2-f_XK^{-1}f_X},\\ 
\bar{f}_t(x) &= m_t(x)+\sigma_t(x)\sqrt{R^2-f_XK^{-1}f_X}, 
\end{aligned}
\end{equation}
where $m_t(x)=f_X^TK^{-1}K_{Xx}$ and $\sigma_t(x)=\sqrt{k(x,x)-K_{xX}K^{-1}K_{Xx}}$. Note that with squared exponential kernel and the sampled points set $X$ to be used in this proof, the invertibility of the matrix $K$ is guaranteed.   
As implied by~\cite[Prop.~1]{maddalena2021deterministic}, 
\begin{equation}\label{eq:bound_truth}
   f(x)\in[\ubar{f}_t(x), \bar{f}_t(x)]. 
\end{equation}
We consider the algorithm that evaluates the grid points $$X=\left\{\left(\frac{k_1}{N}, \frac{k_2}{N},\cdots, \frac{k_d}{N}\right)| k_i\in\{0,1,\cdots, N-1\}\right\}$$ without adaptation, and evaluate the point $\tilde{x}_t$ before termination after $t=N^d$ function evaluations on the grid points, where $\tilde{x}_t$ is given as,
\begin{equation}\label{eq:rep_point_def}
   \tilde{x}_t = \arg\min_{x\in\mathcal{X}} \ubar{f}_t(x).
\end{equation}
Let $x^*$ denote the ground truth optimal solution. We can bound the suboptimality, 
\bse
\label{eq:bound_subopt}
\bee
f(\tilde{x}_t)-\min_{x\in\mathcal{X}} f(x)&\leq&\bar{f}_t(\tilde{x}_t)-f(x^*)\label{inq:bound_r_1}\\
&\leq&\bar{f}_t(\tilde{x}_t)-\ubar{f}_t(x^*)\label{inq:bound_r_2}\\
&\leq&\bar{f}_t(\tilde{x}_t)-\ubar{f}_t(\tilde{x}_t)\label{inq:bound_r_3}\\
&=&2\sigma_t(\tilde{x}_t)\sqrt{R^2-f_XK^{-1}f_X}\\
&\leq&2R\sigma_t(\tilde{x}_t),
\eee
\ese 
where the inequalities in \eqref{inq:bound_r_1} and \eqref{inq:bound_r_2} follow by \eqref{eq:bound_truth} and the inequality \eqref{inq:bound_r_3} follows by the definition of $\tilde{x}_t$ in \eqref{eq:rep_point_def}. We now try to upper bound $\sigma_t(\tilde{x}_t)$. 
We first introduce a set of Lagrangian interpolation functions,
$$
w_{\alpha,N}(x)=\prod_{i\in[N]}\prod_{j\in[N],j\neq \alpha_i}\frac{x_i-j/N}{\alpha_i/N-j/N}, x\in[0,1]^d, \alpha\in[N]^d.
$$
Let $w_N(x)=\left(w_{\alpha,N}(x)\right)_{\alpha\in[N]^d}$. Let $\beta_N(x)=K^{-1}K_{Xx}\in\mathbb{R}^{N^d}$.
We have,
\bee 
\begin{aligned}
(\sigma_{t}(x))^2=&~k(x,x)-K_{xX}K^{-1}K_{Xx}\\
=&~k(x,x)-2K_{xX}\beta_N(x)+\beta_N(x)^TK\beta_N(x)\\
=&~\min_{y\in\mathbb{R}^{N^d}}(k(x,x)-2K_{xX}y+y^TKy)\\
\leq&~k(x,x)-2K_{xX}w_N(x)+w_N(x)^TKw_N(x)
\end{aligned}
\label{eq:bound_posterior_var}
\eee
Let $k_0(x)$ denote the function $k(0, x)$ and $\hat{k}_0$ its corresponding Fourier transformation. By inverse Fourier transformation, we have,
\bee 
\begin{aligned}
&k(x,x)-2K_{xX}w_N(x)+w_N(x)^TKw_N(x)\\
=&~(2 \pi)^{-d} \int_{\mathbb{R}^{d}} \hat{k}_0(\xi)\left(1-2\sum_{\alpha \in {[N]}^d} w_{\alpha, N}(x) e^{i \xi \cdot\left(x-\frac{\alpha}{N}\right)}\right.\\
&~+\left.\sum_{\alpha \in {[N]}^d, \beta \in {[N]}^d} w_{\alpha, N}(x) e^{i \xi \cdot\left(\frac{\beta}{N}-\frac{\alpha}{N}\right)}w_{\beta, N}(x)\right) \diff \xi \\
=&~(2 \pi)^{-d} \int_{\mathbb{R}^{d}} \hat{k}_0(\xi)\left|1-\sum_{\alpha \in {[N]}^d} w_{\alpha, N}(x) e^{i \xi \cdot\left(x-\frac{\alpha}{N}\right)}\right|^{2} \diff \xi \\
=&~(2 \pi)^{-d} \int_{\mathbb{R}^{d}} \hat{k}_0(\xi)\left|e^{-i \frac{\xi}{N} \cdot N x}-\sum_{\alpha \in [N]^d} w_{\alpha, N}(x) e^{-i \frac{\xi}{N} \cdot \alpha}\right|^{2} \diff \xi\\
=&~(2 \pi)^{-d} \int_{\mathbb{R}^{d}} \hat{k}_0(\xi)\left|e^{-i \frac{\xi}{N} \cdot N x}-\sum_{\alpha \in [N]^d} w_{\alpha, N}(x) e^{-i \frac{\xi}{N} \cdot \alpha}\right|^{2} \diff \xi\\
=&~(2 \pi)^{-d} \int_{\xi\in[-\frac{N}{2}, \frac{N}{2}]^d} \hat{k}_0(\xi)\left|e^{-i \frac{\xi}{N} \cdot N x}-\sum_{\alpha \in [N]^d} w_{\alpha, N}(x) e^{-i \frac{\xi}{N} \cdot \alpha}\right|^{2} \diff \xi\\
&~+(2 \pi)^{-d} \int_{\xi\not\in[-\frac{N}{2}, \frac{N}{2}]^d} \hat{k}_0(\xi)\left|e^{-i \frac{\xi}{N} \cdot N x}-\sum_{\alpha \in [N]^d} w_{\alpha, N}(x) e^{-i \frac{\xi}{N} \cdot \alpha}\right|^{2} \diff \xi
\end{aligned}
\label{eq:quad_kernel}
\eee
To proceed, we need to use the Lem.~\ref{lem:bound_fourier}~\cite{zhou2002covering}.
\begin{lemma}[Lemma 4.1, \cite{zhou2002covering}]
Let $x \in[0,1]^{d}$ and $N \in \mathbb{N}$. Then
$$
\sum_{\alpha \in [N]^d}\left|w_{\alpha, N}(x)\right| \leq\left(N 2^{N}\right)^{d}
$$
and for $\theta \in\left[-\frac{1}{2}, \frac{1}{2}\right]^{d}$, there holds
$$
\left|e^{-i \theta \cdot N x}-\sum_{\alpha \in {[N]^d}} w_{\alpha, N}(x) e^{-i \theta \cdot \alpha}\right| \leq d\left(1+\frac{1}{2^{N}}\right)^{d-1}\left(\max _{1 \leq j\leq d}\left|\theta_{j}\right|\right)^{N}.
$$
\label{lem:bound_fourier}
\end{lemma}
We apply the bounds in Lem.~$\ref{lem:bound_fourier}$ to Eq.~\eqref{eq:quad_kernel} and have,
\bee 
\begin{aligned}
&k(x,x)-2K_{xX}w_N(x)+w_N(x)^TKw_N(x)\\
=&~(2 \pi)^{-d} \int_{\xi\in[-\frac{N}{2}, \frac{N}{2}]^d} \hat{k}_0(\xi)\left|e^{-i \frac{\xi}{N} \cdot N x}-\sum_{\alpha \in [N]^d} w_{\alpha, N}(x) e^{-i \frac{\xi}{N} \cdot \alpha}\right|^{2} \diff \xi\\
&~+(2 \pi)^{-d} \int_{\xi\not\in[-\frac{N}{2}, \frac{N}{2}]^d} \hat{k}_0(\xi)\left|e^{-i \frac{\xi}{N} \cdot N x}-\sum_{\alpha \in [N]^d} w_{\alpha, N}(x) e^{-i \frac{\xi}{N} \cdot \alpha}\right|^{2} \diff \xi\\
\leq&~d\left(1+\frac{1}{2^{N}}\right)^{d-1} \max _{1 \leq j \leq d}\left\{(2 \pi)^{-d} \int_{\xi \in\left[-\frac{N}{2}, \frac{N}{2}\right]^{d}} \hat{k}_0(\xi)\left(\frac{\left|\xi_{j}\right|}{N}\right)^{N} \diff \xi\right\} \\
&+\frac{\left(1+\left(N 2^{N}\right)^{d}\right)^{2}}{(2 \pi)^{d}} \int_{\xi \notin\left[-\frac{N}{2}, \frac{N}{2}\right]^{d}} \hat{k}_0(\xi) \diff \xi.
\end{aligned}
\label{eq:quad_kernel}
\eee
We know that $\hat{k}_0(\xi)=(\sigma \sqrt{\pi})^{d} e^{-\frac{\sigma^{2}|\xi|^{2}}{4}}$.
Similar to the analysis in the proof of Example 4 of~\cite{zhou2002covering}, we first try to bound the first term in the upper bound derived in Eq.~\eqref{eq:quad_kernel}.
\bee
\begin{aligned}
&(2 \pi)^{-d} \int_{\xi \in\left[-\frac{N}{2}, \frac{N}{2}\right]^{d}}(\sigma \sqrt{\pi})^{d} e^{-\frac{\sigma^{2}|\xi|^{2}}{4}}\left(\frac{\left|\xi_{j}\right|}{N}\right)^{N} \diff \xi \\
=&\frac{\sigma\sqrt{\pi}}{2\pi} \int_{-N/2}^{N/2} e^{-\frac{\sigma^{2}\xi_j^{2}}{4}}\left(\frac{\left|\xi_{j}\right|}{N}\right)^{N} \left(\prod_{k\neq j}\int_{-N/2}^{N/2}\frac{\sigma \sqrt{\pi}}{2\pi} e^{-\frac{\sigma^{2}\xi_k^{2}}{4}}\diff\xi_k\right)\diff \xi_j \\
\leq&~\frac{\sigma \sqrt{\pi}}{2 \pi} \int_{-N / 2}^{N / 2} e^{-\frac{\sigma^{2}\left|\xi_{j}\right|^{2}}{4}}\left(\frac{\left|\xi_{j}\right|}{N}\right)^{N} \diff \xi_{j} \\
\leq&~\frac{2}{\sqrt{\pi}}\left(\frac{2}{\sigma N}\right)^{N} \Gamma\left(\frac{N+1}{2}\right),
\end{aligned}\label{eq:bound_sigma_1st}
\eee
where $\Gamma(\cdot)$ is the Gamma function. The first inequality in~\eqref{eq:bound_sigma_1st} follows by that $$\int_{-N/2}^{N/2}\frac{\sigma \sqrt{\pi}}{2\pi} e^{-\frac{\sigma^{2}\xi_k^{2}}{4}}\diff\xi_k\leq\int_{-\infty}^{+\infty}\frac{\sigma \sqrt{\pi}}{2\pi} e^{-\frac{\sigma^{2}\xi_k^{2}}{4}}\diff\xi_k=1$$ and the second inequality in~\eqref{eq:bound_sigma_1st} follows by that 
$$
\begin{aligned}
\int_{-N / 2}^{N / 2} e^{-\frac{\sigma^{2}\left|\xi_{j}\right|^{2}}{4}}\left(\frac{\left|\xi_{j}\right|}{N}\right)^{N} \diff \xi_{j}=&2\int_{0}^{N / 2} e^{-\frac{\sigma^{2}\left|\xi_{j}\right|^{2}}{4}}\left(\frac{\left|\xi_{j}\right|}{N}\right)^{N} \diff \xi_{j}\\
\leq&2\int_{0}^{+\infty} e^{-\frac{\sigma^{2}\left|\xi_{j}\right|^{2}}{4}}\left(\frac{\left|\xi_{j}\right|}{N}\right)^{N} \diff \xi_{j}
\end{aligned}
$$ and the definition of Gamma function.
Applying Stirling's formula yields
\bee
\begin{aligned}
&(2 \pi)^{-d} \int_{\xi \in\left[-\frac{N}{2}, \frac{N}{2}\right]^{d}}(\sigma \sqrt{\pi})^{d} e^{-\frac{\sigma^{2}|\xi|^{2}}{4}}\left(\frac{\left|\xi_{j}\right|}{N}\right)^{N} \diff \xi \\
\leq&~\frac{2}{\sqrt{\pi}}\left(\frac{2}{\sigma N}\right)^{N} \Gamma\left(\frac{N+1}{2}\right)\\
\leq&~2\left(\frac{2}{\sigma N}\right)^{N}\left(\frac{N+1}{2 e}\right)^{\frac{N+1}{2}} \frac{1}{\sqrt{N+1}} e^{\frac{1}{6(N+1)}} \\
=&~\sqrt{\frac{2}{e}}e^{\frac{1}{6(N+1)}}\left(\frac{2\sqrt{\frac{N+1}{2}}}{\sigma \sqrt{e}N}\right)^N\\
\leq&~\sqrt{2e}\left(\frac{{2}}{\sigma \sqrt{e N}}\right)^{N}, 
\end{aligned}
\label{eq:upper_bound_first_term}
\eee
where the second inequality in~\eqref{eq:upper_bound_first_term} follows by the Stirling's formula that $$\Gamma(u)\leq \sqrt{2\pi}u^{u-\frac{1}{2}}e^{-u}e^{\frac{1}{12u}}, u>0$$ and the last inequality follows by $e^{\frac{1}{6(N+1)}}\leq e$ and $\frac{N+1}{2}\leq N$.
We are now to bound the second term in Eq.~\eqref{eq:quad_kernel}.
\bee
\begin{aligned}
&(2 \pi)^{-d} \int_{\xi \notin\left[-\frac{N}{2}, \frac{N}{2}\right]^{d}}(\sigma \sqrt{\pi})^{d} e^{-\frac{\sigma^{2}|\xi|^{2}}{4}} \diff \xi \\
=&~\left(\frac{\sigma \sqrt{\pi}}{2 \pi} \int_{\xi_{j} \notin[-\frac{N}{2}, \frac{N}{2}]} e^{-\frac{\sigma^{2}\left|\xi_{j}\right|^{2}}{4}} \diff \xi_{j} \right)^d\\
=&~\left(\frac{\sigma \sqrt{\pi}}{\pi} \int_{N / 2}^{+\infty} e^{-\frac{\sigma^{2}}{4}\left(t^{2}-t / 2\right)} e^{-\frac{\sigma^{2}}{4} \cdot \frac{t}{2}} \diff t\right)^d\\
\leq&~\left(\frac{\sigma \sqrt{\pi}}{\pi} \int_{N / 2}^{+\infty} e^{-\frac{\sigma^{2}}{4}\left(\left(N/2\right)^{2}-\left(N/2\right) / 2\right)} e^{-\frac{\sigma^{2}}{4} \cdot \frac{t}{2}} \diff t\right)^d\\
=&~\left(\frac{\sigma}{\sqrt{\pi}} e^{-\frac{\sigma^{2} N(N-1)}{16}} \frac{8}{\sigma^{2}} e^{-\frac{\sigma^{2} N}{16}}\right)^d\\
=&~\left(\frac{8}{\sigma \sqrt{\pi}}\right)^d e^{-\frac{\sigma^{2}}{16}dN^{2}}.
\end{aligned}
\label{eq:upper_bound_second_term}
\eee
Combining \eqref{eq:quad_kernel}, \eqref{eq:upper_bound_first_term} and \eqref{eq:upper_bound_second_term} yields
\bee
\begin{aligned}
&k(x,x)-2K_{xX}w_N(x)+w_N(x)^TKw_N(x)\\
\leq&~\sqrt{2e}d\left(1+\frac{1}{2^{N}}\right)^{d-1} \left(\frac{{2}}{\sigma \sqrt{e N}}\right)^{N}+\left(1+\left(N 2^{N}\right)^{d}\right)^{2}\left(\frac{8}{\sigma \sqrt{\pi}}\right)^d e^{-\frac{\sigma^{2}}{16}dN^{2}}\\
\leq&~\sqrt{2e}d2^{d-1}\left(\frac{{2}}{\sigma \sqrt{e N}}\right)^{N}+4\left(N 2^{N}\right)^{2d}\left(\frac{8}{\sigma \sqrt{\pi}}\right)^d e^{-\frac{\sigma^{2}}{16}dN^{2}}\\
=&~\sqrt{2e}d2^{d-1}\left(\frac{{2}}{\sigma \sqrt{e N}}\right)^{N}+4\left(\frac{8}{\sigma \sqrt{\pi}}\right)^d e^{-\frac{\sigma^{2}}{16}dN^{2}+2d(N\log2+\log N)}\\
\leq&~\sqrt{2e}d2^{d-1}\left(\frac{{2}}{\sigma \sqrt{e N}}\right)^{N}+4\left(\frac{8}{\sigma \sqrt{\pi}}\right)^d e^{-\frac{\sigma^{2}}{16}dN^{2}+2d(\log2+1) N},
\end{aligned}
\eee
where the first inequality follows by combining \eqref{eq:quad_kernel}, \eqref{eq:upper_bound_first_term} and \eqref{eq:upper_bound_second_term}, the second inequality follows by that $1+\frac{1}{2^N}\leq2$ and $1+(N2^N)^d\leq2(N2^N)^d$, and the last inequality follows by that $\log N\leq N$. Let $N\geq\max\{\frac{32(\log2+2)}{\sigma^2}, \frac{4e^{2d-1}}{\sigma^2}\}$, we have,
\bee
\begin{aligned}
&k(x,x)-2K_{xX}w_N(x)+w_N(x)^TKw_N(x)\\
\leq&~\sqrt{2e}d2^{d-1}\left(\frac{{2}}{\sigma \sqrt{e N}}\right)^{N}+4\left(\frac{8}{\sigma \sqrt{\pi}}\right)^d e^{-\frac{\sigma^{2}}{16}dN^{2}+2d(\log2+1) N}\\
\leq&~\sqrt{2e}d2^{d-1}\left(e^{-d}\right)^{N}+4\left(\frac{8}{\sigma \sqrt{\pi}}\right)^d e^{-dN}\\
=&~\left(\sqrt{2e}d2^{d-1}+4\left(\frac{8}{\sigma \sqrt{\pi}}\right)^d\right)e^{-dN}.
\end{aligned}
\label{eq:bound_the_quad}
\eee
Combining \eqref{eq:bound_posterior_var}, \eqref{eq:bound_the_quad} and that $N^d=t$, we have, 
\begin{equation}\label{ineq:bound_sigma_to_use}
    \sigma_t(x)\leq \left(\sqrt{2e}d2^{d-1}+4\left(\frac{8}{\sigma \sqrt{\pi}}\right)^d\right)^{\frac{1}{2}}e^{-\frac{1}{2}dt^{1/d}}, \forall x\in[0,1]^d.
\end{equation}
Combining that $f(\tilde{x}_t)-\min_{x\in\mathcal{X}}f(x)\leq 2R\sigma_t(\tilde{x}_t)$ in \eqref{eq:bound_subopt} and \eqref{ineq:bound_sigma_to_use}, we have, 
\begin{equation}
    f(\tilde{x}_t)-\min_{x\in\mathcal{X}}f(x)\leq2R\left(\sqrt{2e}d2^{d-1}+4\left(\frac{8}{\sigma \sqrt{\pi}}\right)^d\right)^{\frac{1}{2}}e^{-\frac{1}{2}dt^{1/d}}.
\end{equation}
Setting the right hand side to be smaller than $\epsilon$, we observe that the number of steps $t$ only needs to be $\mathcal{O}\left(\left(\log\frac{R}{\epsilon}\right)^d\right)$. This completes the proof. 
\qed

\section{Proof of Theorem~\ref{thm:mat_lb_ub}}
\thmmatlbub*
\Proof
By Lem.~3 in~\cite{bull2011convergence}, the RKHS on $[0, 1]^d$ is equivalent to Sobolev Hilbert space $H^{\nu+\frac{d}{2}}((0,1)^d)$. Implied by Thm.~1, Sec.~3.3.3~\cite{edmunds1996function}, the covering number of the function space ball in $H^{\nu+\frac{d}{2}}((0, 1)^d)$ is lower bounded by $\Omega\left(\left(\frac{R}{\epsilon}\right)^{\frac{d}{\nu+d/2}}\right)$. Therefore, 
\begin{equation}
    \log\mathcal{N}(S\rX, 4\epsilon, \norm{\cdot}_\infty)=\Omega\left(\left(\frac{R}{\epsilon}\right)^{\frac{d}{\nu+d/2}}\right). 
\end{equation}
We then apply Thm.~\ref{thm:lower_bound} such that we can get the lower bound
\begin{equation}
  T=\Omega\left(\left(\frac{R}{\epsilon}\right)^{\frac{d}{\nu+d/2}}\left(\log\left(\frac{R}{\epsilon}\right)\right)^{-1}\right). 
\end{equation}
The upper bound is implied by Thm.~1 in~\cite{bull2011convergence}. 
\qed


\bibliographystyle{spmpsci_unsrt}. 
\bibliography{ref.bib}


\end{document}